\documentclass[12pt,a4paper,oneside,onecolumn]{article}
\title{Lie algebra of homogeneous operators \\of a vector bundle}
\author{Lecomte P.B.A. and Zihindula Mushengezi E.}
\pdfoutput=1
\usepackage{etex}
\usepackage[]{latexsym}
\usepackage[T1]{fontenc}
\usepackage[latin1]{inputenc}
\usepackage[english]{babel}
\usepackage{graphics}
\usepackage[]{amssymb}
\usepackage{multicol}
\usepackage{amsmath}
\usepackage{graphicx}
\usepackage{amsfonts}
\usepackage[all]{xy}
\usepackage{pifont}
\usepackage{microtype}
\usepackage{answers}
\usepackage{tabularx}
\usepackage{float}
\usepackage{color}

\usepackage{listings}
\definecolor{darkWhite}{rgb}{0.94,0.94,0.94}
\usepackage{palatino}

\newcount\exenum     \newdimen\numindent
\newcount\qqnum     \newdimen\qqmarge
\newcount\sqnum     \newdimen\sqmarge

\outer\def\bye{%
\vskip 0pt\@endmulticol\@endgroup              
\ifanswer \let\next=\exobye@                   
\else     \let\next=\@exobye                   
\fi\next}

\usepackage[Glenn]{fncychap}

\usepackage{pstricks}
\usepackage{pgf,tikz}
\usepackage{pstricks-add}
\usepackage{pst-plot}
\usetikzlibrary{arrows}
\usepackage{pgfplots}
\usepackage{pgfkeys}

\definecolor{qqqqqq}{rgb}{0,0,0}

\definecolor{xdxdff}{rgb}{0.49,0.49,1}
\definecolor{qqwuqq}{rgb}{0,0.39,0}
\definecolor{qqqqff}{rgb}{0,0,1}
\definecolor{ttttff}{rgb}{0.2,0.2,1}
\definecolor{uququq}{rgb}{0.25,0.25,0.25}

\renewcommand{\epsilon}{ \varepsilon}

\newcommand{\euro}{\texteuro{}}

\newcommand{\nit}{{\rm I\!N}}

\newcommand{\pre}{{\bf Proof.\ }}

\newcommand{\cl }{\mathcal }

\newtheorem{theo}{Theorem}[section]
\newtheorem{prop}[theo]{Proposition}
\newtheorem{lem}[theo]{Lemma}
\newtheorem{cor}[theo]{Corollary}

\newcounter{exercice}

\definecolor{fbase}{rgb}{0.8,0.8,1}
\definecolor{fgris}{gray}{0.6}
\definecolor{frouge}{HTML}{DC143C}
\definecolor{fvert}{rgb}{0.6,1,0.6}
\definecolor{fbleu}{rgb}{0.4,0.4,1}
\definecolor{fjaune}{HTML}{DCDC14}
{\endMakeFramed}
\makeatletter
\DeclareRobustCommand\sfrac[1]{\@ifnextchar/{\@sfrac{#1}}%
                                            {\@sfrac{#1}/}}
\def\@sfrac#1/#2{\leavevmode\kern.1em\raise.5ex
         \hbox{$\m@th{\fontsize\sf@size\z@\selectfont#1}$}
         \kern-.1em/\kern-.15em\lower.55ex
          \hbox{$\m@th{\fontsize\sf@size\z@\selectfont#2}$}}

\DeclareRobustCommand{\Efrac}[2]{{\displaystyle\begingroup
\raise2ex\hbox{$\m@th{#1}$}\endgroup\@@over \lower1ex
\hbox{$\m@th{#2}$}}}
\makeatother

\numberwithin{equation}{section}

\newtheorem{Exc}{Exercice}
\Newassociation{correction}{Soln}{mycor}
\Newassociation{indication}{Indi}{myind}

\def\exo#1{\futurelet\testchar\MaybeOptArgmyexoo}
\def\MaybeOptArgmyexoo{\ifx[\testchar \let\next\OptArgmyexoo
                        \else \let\next\NoOptArgmyexoo \fi \next}
\def\OptArgmyexoo[#1]{\begin{Exc}[#1]\normalfont}
\def\NoOptArgmyexoo{\begin{Exc}\normalfont}

\newcommand{\finexo}{\end{Exc}}
\newcommand{\flag}[1]{}

\tikzset{
xmin/.store in=\xmin, xmin/.default=-3, xmin=-3,
xmax/.store in=\xmax, xmax/.default=3, xmax=3,
ymin/.store in=\ymin, ymin/.default=-3, ymin=-3,
ymax/.store in=\ymax, ymax/.default=3, ymax=3,
}

\newcommand{\entete}[1]

\begin{document}
 \Opensolutionfile{mycor}[ficcorex]
 \Opensolutionfile{myind}[ficind]
 \entete{\'Enoncés}


\maketitle

\begin{abstract}
We prove that for a vector bundle $ E \to M$, the Lie algebra $\cl D_{\cl E}(E)$ generated by all differential operators on $E$ which are eigenvectors of $L_{\cl E},$ the Lie derivative in the direction of the Euler vector field of $E,$ and the Lie algebra $\cl D_G(E)$ obtained by Grothendieck construction over the $\mathbb{R}-$algebra $\cl A(E):= {\rm Pol}(E)$ of fiberwise polynomial functions, coincide up an isomorphism.

This allows us to compute all the derivations of the $\mathbb{R}-$algebra $\cl A(E)$ and to obtain an explicit description of the Lie algebra of zero-weight derivations of $\cl A(E).$
\end{abstract}

\section{Definitions}
Let $ \pi: E \rightarrow M $ be a vector bundle, $ L_{\cl E} $ be the Lie derivative in the direction of the Euler vector field and $\cl D (E, M)$ be the Lie algebra of all differential operators $ D: {\rm C} ^ {\infty} (E) \rightarrow {\rm C} ^{\infty} (E),$ where ${\rm C} ^{\infty} (E)$ is the algebra of smooth functions on $E.$\\
A differential operator is said \textit{homogeneous} if it is the sum of differential operators which are eigenvectors with constant eigenvalues of $ L_{\cl E}.$\\
Let now consider the Lie sub-algebra $ \cl D_ {\cl E} (E)$ of $\cl D (E, M) $ of  homogeneous differential operators.\\
Denoting by $ \cl D^k(E, M) $ the space of all differential operators of $ k-$th order,  we have
\[
  \cl D_ {\cl E} (E) = \bigcup_ {k \geq 0} \cl D_{\cl E}^k (E)
\]
where $\cl D^k_{\cl E}(E)$ is the space generated by 
\[\{T\in \cl D^k(E,M)| \exists \lambda\in \mathbb{Z}: L_{\cl E}T=\lambda T\}\]
 The associative algebra $\cl A(E)$ of functions which are fiberwise polynomial coincides with $\cl D_{\cl E}^0(E).$
One has
\[
   \cl A(E)=\bigoplus_{\lambda\in\mathbb{N}}\cl A^{\lambda}(E);
\]
with, for all $u\in\cl A^{\lambda}(E),$ $L_{\cl E}(u)=\lambda u.$\\

The local description of the element of $\cl D_{\cl E}(E)$ is given by the following result taken from \cite{LecLeuZih}.
\begin{prop}\label{local description}
A linear endomorphism $T : \mathcal A(E) \to \mathcal A(E)$ is an element of $\mathcal D_{\mathcal E}(E)$ if and only if in local coordinates $((x^i)_{1\leq i\leq m},(\xi_j)_{1\leq j\leq n})$ corresponding to a local trivialization of $E$, it reads
\[
T=T^{\alpha,\beta}\partial_{\alpha}\overline{\partial}_{\beta} ,
\]
where $\alpha$ and $\beta$ are multi-indices, $\partial_i=\frac{\partial}{\partial x^i}$, $\overline{\partial}_j=\frac{\partial}{\partial \xi_j}$ and all $T^{\alpha,\beta}$ are polynomials in $\xi_1, \dots, \xi_n$ with coefficients in ${\rm C}^{\infty}(M)$.
\end{prop}

We denote by $\cl D_{G}(E)$ the Lie algebra obtained by Grothendieck construction on the associative algebra as follows.
\[
  \cl D_G(E)=\bigcup_{k\geq 0}\cl D_G^k(E),
\] 
where $\cl D_G^0(E)=\cl A(E)$ and 
\[
\cl D^k_G(E)=\{T\in End(\cl A(E)):\forall f\in\cl A(E),[T,f]\in\cl D_G^k(E)\}
\] 
Our aim is to show that 
\[
\cl D_{\cl E}(E)=\cl D_G(E).
\]
 
 \section{The topological space $\cl A(E)$}
 
All open subset $V\subset \mathbb{R}^p$ admits a fundamental sequence of compact  subsets, ie, an increasing sequence $(K_m)$ of compacts in $V$ such that $\cup_m K_m= V$ and $K_m\subset int(K_m)$; where $int(X)$ means the interior of $X.$\\
Let $E$ be a manifold and consider the associative algebra ${\rm C}^\infty(E)$ of all smooth functions. Consider an at most countable atlas $(V_\alpha,\varphi_\alpha)$ of $E$.\\
 Let $(K_{m,\alpha})$ be a fundamental sequence of compact subsets of $\varphi_\alpha(V_\alpha).$\\
 
 For all $(s,m)\in \mathbb{N}\times\mathbb{N}_0$ and $f\in {\rm C}^\infty(E),$ we put 
 \[
  p_{s,m,\alpha}(f)= \sup_{x\in K_{m,\alpha},|\lambda|\leq s}
                             |\partial^\lambda(f\circ\varphi_\alpha^{-1})|,
 \]
with $\lambda $ a multi-index. \\
These $ p_{s,m,\alpha}$ as defined above are semi-norms on ${\rm C}^\infty(E)$. 
They provide with the space ${\rm C}^\infty(E)$ an Hausdorff, locally convex and complete topological space structure, i.e, a FRECHET space structure.
This topology has the following property.\\
A sequence $(f_k)$ of functions in ${\rm C}^\infty(E)$ converges to zero if and only if for all chart $(V,\varphi)$ of $E,$ for all compact subset $K$ of $\varphi(V)$ and all multi-index $\lambda,$ the sequence of all restrictions of $\partial^\lambda(f\circ\varphi_\alpha^{-1})$ to $K$ converge uniformly to zero.\\
Since topology defined on ${\rm C}^\infty(E)$ is associated with a countable family of semi-norms, then ${\rm C}^\infty(E)$ is metrizable. \\
 In the following lines, the space $\cl A(E)$ is provided with the topology induced by that of ${\rm C}^\infty(E)$.\\
As a result, a function $\Phi: \cl A(E)\rightarrow \cl A(E)$ is continuous if for every sequence $(P_n)$ in $\cl A(E)$, $\Phi(P_n)\rightarrow \Phi(P)$ whenever $P_n\rightarrow P.$\\
 
Generally, if a topological vector space $ V $ is a direct sum of two vector subspaces $ V_1 $ and $ V_2, $ this is not enough to say that $ V $ is also their topological sum.\\  
Let $ \pi: E \rightarrow M $ be a vector bundle. \\ 
\begin{prop}
For all $r\in \mathbb{N},$ the linear application 
 \[
 pr_r: \cl A(E)\rightarrow \cl A^r(E): u=u_0+u_1+...+u_s\mapsto u_r
 \]
 is continuous.
\end{prop}
\pre
Let $(P_n)$ be a sequence that converges to zero in $\cl A(E).$  
Let $K$ be a compact in $\varphi(V),$ where $(V,\varphi)$ is an adapted chart of the vector bundle $E.$ 
Consider a fundamental sequence $(K_m)$ of compacts in $\varphi(V)$ with $0\in K_m$ . 

There exists $m\in\mathbb{N}$ such that $K\subset K_m.$
Locally, one can write 
\[
  P_n=\sum_{i=0}^{d_n}P_n^i,
\]
where $(P_n^i\circ\varphi^{-1})(x,\xi)=\sum_{|\alpha|=i}A_n^\alpha(x)\xi^\alpha$ and $d_n\in\mathbb{N}.$ 
We will show that the sequence $(P_n^r)$ also converges to zero for all $r\in\mathbb{N}.$\\
Let $\epsilon >0.$ There exists $N\in\mathbb{N}$ such that 
\begin{equation}\label{(*)}
 \forall n\geq N, \sup_{(x,\xi)\in K_m,|\lambda|+|\mu|\leq s}
 |\partial^\lambda\overline{\partial}^\mu (P_n\circ\varphi^{-1})(x,\xi)|<\epsilon 
\end{equation}
Observe that for any $r\in\mathbb{N}_0,$ we have
\begin{eqnarray*}
   p_{s,m,\lambda,\mu}(P_n^r\circ\psi^{-1}(x,\xi)) & = &  \sup_{(x,\xi)\in K_m,|\lambda|+|\mu|\leq s}
| \sum_{|\alpha|= r}\partial^\lambda( A_n^
    \alpha)(x)\overline{\partial}^\mu \xi^\alpha| \\
                             & \leq & \sup_{(x,\xi)\in K_m,|\lambda|+|\mu|\leq s}
\sum_{|\alpha|= r}| \partial^\lambda( A_n^
    \alpha)(x)||C_{\alpha,\mu}| 
\end{eqnarray*}
with $C_{\mu,\alpha}$ the maximum of the function $(x,\xi)\mapsto \overline{\partial}^\mu \xi^\alpha$ in $K_m.$\\
We therefore have to show that for any given multi-index $\beta$ such that $|\beta|=r,$  the sequence $(A_n^{\beta})$ converges to zero. 
From the relation  (\ref{(*)}),  we deduce that from a certain rank, we can write
  \begin{eqnarray*}
  \varepsilon\beta! & > & \sup_{(x,\xi)\in K_m,|\lambda|+|\mu|\leq s+r}
                   |\partial^\lambda\overline{\partial}^\mu (P_n\circ\psi^{-1})(x,\xi)|\\  
                   & = & \sup_{(x,\xi)\in K_m,|\lambda|+|\mu|\leq s+ r}|\sum_{| 
                   \alpha|\leq d_n}\partial^\lambda A_n^\alpha(x)\overline{\partial}^\mu \xi^\alpha|\\
                   & \geq & \sup_{(x,0)\in K_m,|\lambda|+|\mu|\leq s+ r,\mu=\beta} 
                 |\sum_{|\alpha|\leq d_n}\partial^\lambda A_n^\alpha(x)\overline{\partial}^\mu \xi^\alpha|\\
                   & = &  \beta!\sup_{(x,\xi)\in K_m,|\lambda|\leq s}|\partial^\lambda A^\beta_n(x)| 
 \end{eqnarray*}
 This allows  then to conclude that $(P_n^r)\to 0.$  \hfill $\blacksquare$
\section{Algebras identification : $\cl D_{\cl E}(E)\cong\cl D_G(E)$}
We first start by showing that given a vector bundle $ \pi: E \to M,$ a differential operator of $ \cl D (E)$ is entirely determined by its values on the fiberwise polynomial functions.

Let us state the following preliminary result which can be justified as in \cite{Lec4} p.7; the main ingredient of the proof being the Taylor development.

\begin{lem}\label{lem1}
 If $ u\in\mathcal{ A}(E)$ is such as  $ j_{_a}^l u=0,$ then we can write, in the neighborhood of $a,$
    \begin{equation*}\label{(i)}
      u=\sum_{i=1}^N u_{i_0}\cdots u_{i_l} 
    \end{equation*}
  with $u_{i_j}\in\mathcal{ A}(E)$ and $u_{i_j}(a)=0,$ 
                                $\forall (i,j)\in [1,N]\times[0,l].$  
\end{lem}
The following proposition will allow us to extend each element of $ \cl D_G(E) $ into a differential operator of $\cl D(E).$
 \begin{prop}\label{prop D vers D chapeau}
  For all endomorphism $D:\cl A(E)\to \cl A(E)$ of $\mathbb{R}-$vector spaces such that 
  \[
   j_a^l(u)=0 \Rightarrow D(u)(a)=0,\quad\forall u\in\cl A(E),
  \]
 there exists $\widehat{D}\in\cl D(E)$ such that 
   \[
     \widehat{D}(v)=D(v), \quad\forall v\in\cl A(E).
   \]
 \end{prop}
\pre
This statement results from the fact that for all $f\in {\rm C}^\infty(E),a\in E $ and for all integer $k\in\mathbb{N}, $ there exist $u\in\cl A(E)$ such that 
   \[
     j_{_{a}}^k(f)=j_{_{a}}^k(u).
   \]
Let then $D\in\cl D_G^k(E).$ One can therefore set, for all  $f\in {\rm C}^\infty(E),$ 
 \[
  \widehat{D}(f)_{_{a}}=D(u)_{_{a}}, 
 \] 
the polynomial function $u$ having the same $k$-order jet as $f$ on $a\in E.$  \hfill $\blacksquare$ 

\begin{cor}\label{lem2}
  For all $D\in\cl D^l_{_G}(E)$, $(l\geq0),$ there is a unique differential operator   
  $\widehat{D}\in \cl D^l(E,M)$ such that  
    \[
      \widehat{D}(u)=D(u),\,\forall u\in\cl A(E).
    \]   
\end{cor}
\pre
Let $D\in\cl D^l_{_G}(E).$ Consider a function  $u\in\cl A(E)$ such that $j_{_{a}}^l(u)=0.$ 
Let us show that we have $D(u)_{_{a}}=0.$ We proceed by induction on $l.$ Indeed, this statement being true for $l=0,$ suppose, by induction hypothesis that it is for $k<l.$ So, when $k=l,$ with the notations of the previous Lemma \ref{lem1}, we have that
  \[
    D(u)= \sum_{i=1}^N u_{i_0}D(u_{i_1}\dots u_{i_l})+ 
  \sum_{i=1}^N \underbrace{[D,u_{i_0}]}_{\in\cl D_{_G}^{l-1}(E)}(u_{i_1}\dots u_{i_l})
  \]  
vanishes on $a$ and then, the desired result follows. \hfill $\blacksquare$\\

At this point we have already established that the elements of $ \cl D_G(E),$ like those of $ \cl D_{\cl E} (E),$ can be seen as restrictions of the differential operators  of $ \cl D (E) $ on the $\mathbb{R}-$algebra $\cl A(E).$\\
Thus, the elements of $ \cl D_G (E) $ locally decompose into expressions comprising polynomial functions along the fibers of $ E $. \\
But nothing tells us that these polynomials are of bounded degree independently of charts; as it turned out to be the case with the elements of $\cl D_{\cl E}(E)\footnote{It shown in \cite{LecLeuZih}}$.\\
The following lemma will allow us to prove a result which states that this is also true for the elements of $ \cl D_G (E). $  
  
\begin{lem}\label{lem3}
 Let $D\in\cl D^l_{_G}(E).$  There is $N_s\in\mathbb{N}$ such that 
    \[
      \cl A^s(E)=\cl F_{N_s}(E),
    \]
 where, for $N\in\mathbb{N},$ by definition, one has
    \[
     \cl F_{N}(E)=\{u\in\cl A^s(E):r\geq N\Rightarrow pr_r\widehat{D}(u)=0)\}\cdot  
    \] 
\end{lem}
\pre

Note that 
 \[
 \cl A^s(E)=\{u\in {\rm C}^\infty(E): L_{\cl E}u=su\} 
 \]
is closed in ${\rm C}^\infty(E);$ it is therefore a Baire space. For $N\in\mathbb{N},$ 
it is the same for
\[
  \cl A^s_{_N}(E)=\{u\in\cl A^s(E): r\geq N\Rightarrow pr_{_r}\widehat{D}(u)=0\}.
\] 
This last set is an intersection of closed subsets, due to the continuity of the differential operators and that of $ pr_r,$ by virtue of the previous Proposition \ref{lem2}. (Relative to the topology for which $ {\rm C}^\infty (E) $ is a Fréchet space.)\\
Observe that we have
  \[
    \cl A^s(E)=\bigcup_{_{N\in\nit}}\cl A^s_{_N}(E)\cdot
  \]
Thus, there is $N_s\in\mathbb{N}$ such that $\cl A^s_{_{N_s}}(E)$ is of non-empty interior. However, any open set of a topological vector space containing the origin is absorbing.
So we can write
 \[
   \cl A^s(E)=\cl A^s_{_{N_s}}(E)\cdot
 \]
Indeed, let $v\in Int(\cl A_{N_s}^s(E)).$ There is then an open set $\cl U$ containing the origin such that 
 \[
    v+\cl U\subset Int(\cl A_{N_s}^s(E))\subset\cl A^s(E)\cdot
 \]
Now consider an element $u\in\cl A^s(E).$ As $\cl U$ is absorbent, there is $\kappa>0$ such that $\kappa u\in\cl U.$ Therefore,  
\[
\kappa u+v\in\cl A^s(E)
\]
and thus $u\in\cl A^s(E).$  \hfill $\blacksquare$

\begin{prop}\label{caract}
   
The unique differential operator $\widehat{D}\in\cl D(E)$ associated with $D\in\cl D^l_{_G}(E)$ is written locally
   \begin{equation}\label{(**)}
      \widehat{D}=\sum_{\left|\alpha\right|+\left|\beta\right|\leq l} 
                         u^{\alpha,\beta}\partial_\alpha\overline{\partial}_\beta 
   \end{equation}
 where the $u^{\alpha,\beta}$ are polynomials in $ (y^1, \dots, y^n) $ of a maximum degree bounded independently of charts.
\end{prop}
\pre
Recall first that
   \[
     \partial_\alpha\overline{\partial}_\beta(x^\gamma y^\delta)=\left\{\begin{array}{l}
                                          \frac{\gamma!}{(\gamma-\alpha)!}\frac{\delta!}{(\delta-\beta)!}     
              x^{\gamma-\alpha}y^{\delta-\beta}\mbox{ if } \alpha\leq \gamma 
              \mbox{ and } \beta\leq\delta\\
                         0 \quad\quad\mbox{ if not }.
               \end{array} \right.                                                                 
   \] 
Let's do the proof by induction. By virtue of the previous Lemma\ref{lem3}, for $\left|\alpha\right|+\left|\beta\right|=0,$ there exits $N_0\in\mathbb{N},$ such that
   \[
     u^{00}=\widehat{D}(1)\in\cl A^0\oplus\dots\oplus\cl A^{N_0}(E).
   \]
Assume by induction hypothesis that the proposition is true for $\left|\alpha\right|+\left|\beta\right|<t.$ 

Let $ U $ be a trivialization domain over which $ \widehat{D} $ is of the form (\ref{(**)}). \\
Consider a function $\rho\in\cl A^0(E)$ such that $\rho=\pi_E^*(f), f$ being a function with compact support in $ U, $ zero outside $ U $ and equal to $ 1 $ in an open $V$ of $U.$ \\
For $\left|\gamma\right|+\left|\delta\right|=t,$ we have
  \[
   \widehat{D}(\rho\,x^\gamma y^\delta)\stackrel{V}{=}
   \sum_{\stackrel{|\alpha|+|\beta|\leq l}{\alpha<\gamma \mbox{ or } \beta<\delta}} u^{\alpha,\beta}(x,y)\frac{\gamma!}{(\gamma-\alpha)!}\frac{\delta!}{(\delta-\beta)!}  
             x^{\gamma-\alpha}y^{\delta-\beta}+\gamma!\delta! u^{\gamma,\delta}(x,y)
  \]
Observe that there is $N_{|\delta|}$ such that
   \[
\widehat{D}(\rho\,x^\gamma y^\delta)\in
                                      \cl A^0(E)\oplus\dots\oplus\cl A ^{N_{\left|\delta\right|}}(E). 
   \]                                    
This comes from the previous Lemma \ref{lem3}, and we have that for $\alpha<\gamma\mbox{ or }\beta<\delta,$ the recurrence hypothesis can be applied to  $u^{\alpha,\beta}.$\\
We deduce that $u^{\gamma,\delta}$  is polynomial of bounded degree, independently of chart; and this completes the proof of the proposition. \hfill $\blacksquare$\\

From the previous proposition and from the local characterization of the algebra $\cl D_{\cl E}(E)$  we deduce the following result.

\begin{theo}
  Let $E\to M$ be a vector bundle. 
Quantum Poisson algebras $\cl D_{_G}(E)$ and $\cl D_{\cl E}(E)$ coincide up to isomorphism.
\end{theo}

\section{Derivations of the associative algebra $\cl A(E)$}
Let $ E \to M $ be a vector bundle. In the following lines, $ \cl A (E) $ still designates the associative algebra of polynomial functions along the fibers of the bundle $ E \to M. $
We use the equality $ \cl D_{\cl E} (E) = \cl D_G (E) $ to determine all the derivations of the  $\mathbb{R}-$algebra $ \cl A (E).$ We denote by $Vect(E)$ the space of vector fields of $E.$
  
\begin{prop}
A linear map $D:\cl A(E)\to \cl A(E)$ is a derivation of $\cl A(E)$ if and only if $D$ is the restriction to $\cl A(E)$ of an element of  $\cl D^1_{\cl E}(E)\cap Vect(E).$\\
In other words,
 \[
 Der(\cl A(E))=\cl D^1_{\cl E}(E)\cap Vect(E)|_{\cl A(E)}.
 \]
\end{prop}

\pre
The inclusion $Der(\cl A(E))\supset\cl D^1_{\cl E}(E)\cap Vect(E)$ is obvious.
Now let $D\in Der(\cl A(E)).$ We have, for all $u\in\cl A(E),$ the following equality
\[
 [D,\gamma_u]=\gamma_{D(u)}\cdot
\]
Indeed, for all $ u, v \in \cl A (E), $ we can write
\begin{eqnarray*}
    [D,\gamma_u](v)  & = & D(uv)-uD(v)\\
              & = & D(u)v= \gamma_{D(u)}(v)\cdot
\end{eqnarray*}
Therefore, $D\in\cl D_G^1(E)=\cl D^1_{\cl E}(E).$
Observe that as $\mathbb{R}-$vector spaces, we have the equality
 \[
   \cl D^1_{\cl E}(E)=(\cl D^1_{\cl E}(E)\cap Vect(E))\oplus \cl A(E)\cdot
 \] 
 Let us then set $D=D_c+w$ with $D_c\in\cl D^1_{\cl E}(E)\cap Vect(E)$ and $w\in\cl A(E).$ As $D(1)=D_c(1)=0,$ we conclude that $w=0.$ Therefore, we obtain the following inclusion  
 \[
 Der(\cl A(E))\subset\cl D^1_{\cl E}(E)\cap Vect(E);
 \]
 which completes the demonstration. \hfill $\blacksquare$\\

We propose in the following lines a result which relates the Lie algebra of the infinitesimal automorphisms $ Aut(E) $ of the vector bundle $ E $ and that of zero-weight derivations of the  $\mathbb{R}-$algebra $\cl A(E).$ For the proof, see in \cite{LecLeuZih}.

\begin{prop}\label{der0}
 The algebra of zero-weight derivations of $ \cl A (E) $ is given by
 \[
   Der^0(\cl A(E))=Aut(E)_{|\cl A(E)}\cdot
 \]
\end{prop}

This proposition is used, in \cite{LecLeuZih}, to prove a Lie-algebraic characterization of vector bundles, by virtue of a result taken from \cite{Lec1}.\\

We also know that given two vector bundles $E\to M$ and $F\to N,$ any isomorphism  $\Psi:\cl A(E)\to \cl A(F)$ of associative algebras induces a Lie algebra  isomorphism by
\begin{equation}\label{psi ind der}
\widehat{\Psi}: Der(\cl A(E))\to Der(\cl A(F)): D\mapsto \Psi\circ D\circ \Psi^{-1}.
\end{equation}

Moreover, if $ \Psi $ is graded, the induced isomorphism respects the Lie subalgebras $Der^0(\cl A(E))$ and $Der^0(\cl A(F))$  of zero-weight derivations of    $Der(\cl A(E))$ and $ Der(\cl A(F))$ respectively.\\
We propose in this section another way to determine the zero-weight derivations of $\cl A(E)$ and in doing so, we show that such an isomorphism $\widehat{\Psi},$ given in (\ref{psi ind der}), preserves the Euler vector field.

\begin{theo}\label{der poids nul}
Let $\pi:E\longrightarrow M$ be a vector bundle of rank $n$.\\ 
  (a) The Lie algebra of homogeneous zero-weight derivations of the associative algebra $\cl A(E)$ is given by
    \[
      Der^0(\cl A(E))\cong Vect(M)\oplus gl(E^*),     
    \] 
where $gl(E^*)=\Gamma(Hom(E^*,E^*))$ is the space of all smooth sections of the vector bundle $Hom(E^*,E^*)\to M.$\\
  (b) The center $Z(Der^0(\cl A(E)))$ of this Lie algebra is formed by the real multiples of the Lie derivative in the direction of the Euler vector field.   
\end{theo}
\pre
Let $D\in Der^0(\cl A(E)).$ Since $D$ respects the gradation of $\cl A(E),$ its restriction to $\cl A^0(E)=\{f\circ\pi: f\in{\rm C}^\infty(M)\}$ comes down to the action of a vector field $X\in Vect(M)$ by 
    \begin{equation}\label{der vect}
    \pi^*f\mapsto \pi^* L_X f\cdot
    \end{equation}     
Consider the restriction of $D$ to $\cl A^1(E),$ this last space being identified with $\Gamma(E^*).$ In fact, with $u\in\Gamma(E^*)$ we associate $\widehat{u}\in\cl A^1(E)$ defined by
  \[
    \widehat{u}(a)=u_x(a),
  \]  
for $a\in E_x.$ 
Therefore, the following linear application 
\[
D:\Gamma(E^*)\longrightarrow \Gamma(E^*): u\mapsto D(u):=D(\widehat{u})
\]
is a first-order differential operator.

Indeed, for any $u\in\Gamma(E^*),$ such that $j^1_x u=0$ in an open $U$ of $M$ containing $x,$ we can consider a decomposition
\[
u=\Sigma_i f_iu_i, \quad f_i\in\cl A^0(E),u_i\in\Gamma(E^*)
\]
 where the $f_i$ and the $u_i$ vanish on $x.$ 
Therefore, the differential operator $D$ acting on the sections of the bundle $ E^*$ is written locally
  \[
     D(u)= A(u)+ \Sigma_i A^i(\partial_iu)
  \]
with $A^i,A\in {\rm C}^{\infty}(U,gl(n,\mathbb{R})).$ 
Observe that $\Gamma(E^*)$ is a $\cl A^0(E)-$module; we define then for any $f\in\cl A^0(E),$ 
a zero-order differential operator acting on the sections of the bundle $E^*$ by
 \[
   \gamma_f: \Gamma(E^*)\to\Gamma(E^*): u\mapsto fu .
 \]
We then have on the one hand,
\begin{eqnarray*}
   D(fu)& = & fD(u)+ D(f)u\\   
        & = & fA(u)+f\Sigma_i A^i(\partial_i u) + D(f)u,
 \end{eqnarray*} 
 for any $f\in\cl A^0(E)$ and $u\in\Gamma(E^*)\cong \cl A^1(E).$ 
 And on the other hand,
 \begin{eqnarray*}
   D(fu)& = & A(fu)+ \Sigma_i A^i(\partial_i(fu))\\   
        & = & A(fu)+\Sigma_i A^i(f\partial_i u) + \Sigma_i A^i(u \partial_i f).
 \end{eqnarray*} 
In addition, by using the relation (\ref{der vect}), we have 
 \begin{eqnarray*}
   D(f)u& = & (X.f)u ,
 \end{eqnarray*} 
for a certain vector field $X\in Vect(M).$ \\
And more, for all $f\in\cl A^0(E)$ and $A\in gl(E^*),$ we have that $[A,\gamma_f]=0,$ this last bracket being that of the commutators in the algebra of endomorphisms of the $\mathbb{R}-$vector space $\Gamma(E^*).$\\
We deduce that $ A^i = X^i id $ and we can therefore write 
 \[
   D(u)=\nabla_{X}u+B(u),
 \]
with $X\in Vect(M)$ and $B\in gl(E^*).$ \\
We have assumed given, in what precedes, a connection on the vector bundle $ E \to M, $ and this is still the case in the following lines.\\
Observe that the derivation $L_{\cl E}\in Der^0(A(E))$ is zero in $A^0(E)$ and that it coincides with the identity on $\Gamma(E^*);$ which therefore corresponds to the case $X=0 , A=id.$ \\
Since $\cl A^0(E)$ and $\cl A^1(E)$ generate the all $\mathbb{R}-$algebra $\cl A(E),$ the part $(a)$ of the theorem is thus established.

Let $D_{X,A}\in Z(Der^0(\cl A(E))), D_{Y,B}\in Der^0(\cl A(E))$ and $u\in\Gamma(E^*).$\\ We must have
  \begin{eqnarray*}
   0 & = & [D_{X,A},D_{Y,B}](u) \\ 
     &=&  (R(X,Y)+ \nabla_{[X,Y]})u + (\nabla_X B)(u)-(\nabla_Y A)(u) 
                             +   [A,B](u)
  \end{eqnarray*}                           
 This relation being true for all $ Y $ and all $ B,$ we obtain, by setting $ Y = 0,$ 
   \[
    [A,B]=0 \mbox{ and } \nabla_X B= 0, \qquad\forall B.
   \]
The first equality gives $A=\kappa\,id.$ We deduce from the second that, by setting  
 $B=f\,id,$ with $f\in\cl A^0(E),$ 
 \[ \left( X.f \right) id=0 \] 
 and thus $X=0.$ \hfill $\blacksquare$

\begin{cor}
  Let $E\to M$ and $F\to N$ be two vector bundles. 
 
Any isomorphism of associative algebras $ \Phi: \cl A(E) \to \cl A(F) $ induces an isomorphism of Lie algebras $\widehat{\Psi}: Der(\cl A(E))\to Der(\cl A(F))$ such that
  \[
    \widehat{\Psi}(\cl E_E)=\cl E_F \quad \mbox{ and } \quad\widehat{\Psi}(Der^r(\cl A(E)))= Der^r(\cl A(F))
  \]
with $Der^r(\cl A(E))$ (resp.$ Der^r(\cl A(F)))$ designating the $\mathbb{R}-$vector  space of $r$-weight derivations of $\cl A(E)$ (resp. $\cl A(F))$.
\end{cor}
\pre
   Since $\Psi$ induces a graded isomorphism between $\cl A(E)$ and $\cl A(F),$ the proof of this can be founded in \cite{LecLeuZih},
we will denote them both by $\Psi.$ Therefore, by definition, for all  $D\in Der(\cl A(E)),$ we have 
\[\widehat{\Psi}(D)=\Psi\circ D\circ\Psi^{-1};\] 
and thus, $\Psi$ respecting the gradations, it comes, for all $D\in Der^r(\cl A(E))$ and all $u\in\cl A^s(F),$
   \begin{eqnarray*}
       (\widehat{\Psi}(D))(u) & = & \Psi(D(\Psi^{-1}(u))) \in\cl A^{r+s}(F),
   \end{eqnarray*}
 because $\Psi^{-1}(u)\in\cl A^s(E).$ 
Hence the following inclusion 
\[
\widehat{\Psi}(Der^r(\cl A(E)))\subset Der^r(\cl A(F)).
\]
In addition, as a result, we have
 \[
 \widehat{\Psi}(Z(Der^0(A(E))))=Z(Der^0(A(F))).
 \]
This comes from the previous Theorem \ref{der poids nul} and allows to conclude that there exists $\kappa\in\mathbb{R}\setminus\{0\}$ such that
  \[
  \widehat{\Psi}(\cl E_E)=\kappa \cl E_F .
  \]
Therefore, for any $u\in\cl A^1(F),$ on the one hand we have
   \begin{eqnarray*}
       (\widehat{\Psi}(\cl E_E))(u) & = & \Psi\circ \cl E_E\circ\Psi^{-1}(u) \\
                                    & = & \Psi(\Psi^{-1}(u)) \quad\mbox{ since } \Psi^{-1}(u)\in\cl A^1(E)
   \end{eqnarray*}
 And on the other hand,  
   \begin{eqnarray*}
       (\widehat{\Psi}(\cl E_E))(u) & = & (\kappa \cl E_F)(u) \\
                                    & = & \kappa u
   \end{eqnarray*}
 Therefore, $\kappa=1;$ and we have that $\widehat{\Psi}$ preserves the Euler vector field. \hfill $\blacksquare$

\end{document}